\documentclass[pre,twocolumn,showpacs,superscriptaddress,preprintnumbers,floatfix]{revtex4}

\usepackage{amssymb,amsmath,latexsym,epsf}
\usepackage{dcolumn}
\usepackage{graphicx}
\usepackage{bm} 

\newcommand{\eM}     {$\epsilon$-machine}
\newcommand{\eMs}    {$\epsilon$-machines}

\newcommand{\biinfinity}	{ \stackrel{\leftrightarrow} {s} }

\newcommand{\past}	{ {\stackrel{\leftarrow} {s}} }
\newcommand{\pastprime}	{ {\past}^{\prime}}

\newcommand{\future}	{ \stackrel{\rightarrow}{s} }

\newcommand{\AllPasts}	{ { \stackrel{\leftarrow} {\rm {\bf S}} } }

\newcommand{\CausalState}	{ {\cal S} }

\newcommand{\CausalStateSet}	{ \boldsymbol{\CausalState} }

\newcommand{\Prob}		{ {\rm P}}

\newcommand{\Cmu}		{{C_\mu}}

\newcommand{\sAlg}   {$\sigma$-algebra}

\newcommand{\OSym}   { \Delta }

\begin{document}


\title{Reductions of Hidden Information Sources}

\author{Nihat Ay}
\email{nay@santafe.edu}
\affiliation{Santa Fe Institute, 1399 Hyde Park Road, Santa Fe, NM 87501}
\affiliation{Mathematical Institute, Friedrich-Alexander 
University Erlangen-Nuremberg, Bismarckstr. 1, 91054 Erlangen, Germany}

\author{James P. Crutchfield}
\email{chaos@santafe.edu}
\affiliation{Santa Fe Institute, 1399 Hyde Park Road, Santa Fe, NM 87501}

\date{\today}

\bibliographystyle{unsrt}


\begin{abstract}
In all but special circumstances, measurements of time-dependent processes
reflect internal structures and correlations only indirectly. Building
predictive models of such hidden information sources requires discovering,
in some way, the internal states and mechanisms. Unfortunately, there are
often many possible models that are observationally equivalent. Here we show
that the situation is not as arbitrary as one would think. We show that
generators of hidden stochastic processes can be reduced to a minimal form
and compare this reduced representation to that provided by computational
mechanics---the \eM. On the way to developing deeper, measure-theoretic
foundations for the latter, we introduce a new two-step reduction process.
The first step (internal-event reduction) produces the smallest observationally
equivalent \sAlg\ and the second (internal-state reduction) removes
\sAlg\ components that are redundant for optimal prediction. For several
classes of stochastic dynamical systems these reductions produce
representations that are equivalent to \eMs.
\end{abstract}

\pacs{
  02.50.Ey, 
  02.50.Ga, 
  05.45.Tp, 
  89.70.+c  
  \\
\begin{center}
  Santa Fe Institute Working Paper 04-05-011;
  arxiv.org/abs/cs.XX/0405XXX
\end{center}     
  }
\maketitle


\tableofcontents

\section{Introduction}

Experiment and simulation often produce voluminous amounts of data---data that
the scientist or analyst attempts to understand by building predictive models.
The best models, however, do more than simply predict the data. In the best of
circumstances, models also capture the internal structures, active degrees of
freedom, correlations, and so on that underlie the observations. In this way
modeling enhances understanding and leads to new insights about the forces
that shape our world.

Unfortunately, measurements generally are only indirect indicators of internal
structure. This makes the process of model building difficult and often highly
nonunique. One would hope that there is some principled approach to model
building and inference that would guide us in inferring structural
properties from data. The possibilities for such an
approach are bounded by two extremes: (i) Are there formal constraints that
guide the discovery of good representations? (ii) Can the observations
themselves tell us which representation to use or, perhaps, how to correct
an initially faulty hypothesis?

These days, however, the problem of building useful predictors of hidden
information sources is compounded by the fact that the systems studied are
quite complicated, in the sense of consisting of many components, for example.
Genomic, geophysical, neurobiological, Internet traffic, and World Wide Web
systems easily come to mind as complex in this sense and as particularly
desirable to model. This very practical observation, in turn, argues even
more forcefully for a principled approach to discovering and describing hidden
structure. That is, we now need to understand the process of model building
for such complicated systems well enough to teach machines how to do it.

Here, building on previous work \cite{Crut88a,Crut98d,Shal98a}, we address one
piece in this puzzle---what we call the \emph{Forward Modeling Problem}: Given
a generator of an observed stochastic process, Is there a minimal, optimal
predictor of it? In answering this question positively, we have two goals.
The first, naturally enough, is to articulate the notion of minimal generators
of observed stochastic processes and show that they exist. The second, though,
is to lay more rigorous and broader foundations than currently available for
the \emph{Reverse Modeling Problem}---Given observations, can one reconstruct
the hidden mechanisms?

\subsection{Background}

Reviewing a little background on previous work will help put the current formal
results in perspective and motivate our development. We then comment on closely
related work in which similar questions arise, but which take different
approaches to structural inference. Then, after outlining our approach, the
mathematical development begins.

If we are to build a predictive model of an information source that produces
a time series, the most basic assumption to make is that the source, at each
moment of time, is in some ``state''. Over time, the source transitions from
state to state. As we noted already, though, in the general setting we do not
have access to these states, we only have indirect information about
them---information that we call \emph{measurements}. So the modeling question
reduces to the following. Given that all we have are sequences of observations,
what kind of ``state'' should be formed from them and used for modeling? The
answer is rather straightforward, and seemingly tautological: the ``states''
that we should use are those that are effective for prediction.

This is the starting point for how \emph{computational mechanics}
\cite{Crut88a,Crut98d,Shal98a} builds optimal models. One of
the notable results in computational mechanics, though, is that the
representation, which emerges from focusing on states that are effective for
prediction, captures all of a process's internal causal structure.
In fact, computational mechanics shows that there is a preferred
representation for modeling, which is called an \emph{\eM}.

To start, we consider a time series of observations
$\biinfinity = \ldots, s_{-2}, s_{-1}, s_0, s_1, \dots$, in which the
individual measurements are symbols in a finite alphabet: $s_i \in \OSym$.
An \eM\ consists of states---called \emph{causal states} and denoted
$\CausalStateSet$---and transitions between them. The causal states are
defined as those \emph{sets of histories}
$\past_t = \ldots, s_{t-3}, s_{t-2}, s_{t-1} $
that are equivalent for predicting the \emph{future}
$\future_t = s_t, s_{t+1}, s_{t+2}, \ldots$. That is, two histories---$\past$
and $\pastprime$---are associated with a given causal state, when the
\emph{sets} of possible futures ``look'' the same having seen them. More
precisely, this modeling principle defines an equivalence relation $\sim$
over the set $\AllPasts$ of histories:
\begin{equation}
\past \sim \pastprime ~\mathrm{if~and~only~if}~
  \Prob ( \future | \past ) = \Prob ( \future | \pastprime ) ~,
\end{equation}
where in the conditional distribution equality we mean that each individual
future is given the same probability. The resulting equivalence
classes are the causal states.

From this, one can show that the \eM\ for an information source is the optimal,
minimal, and unique predictor of an information source. In the language of
mathematical statistics, the \eM\ is a minimal sufficient statistic for the
observed stochastic process produced by an information source. More than being
a good predictor that is small, the semigroup determined by the causal states
and transitions captures all of the information source's internal
structure---regularities, symmetries, and so on. And, due to minimality, one
can show that the \emph{statistical complexity} $\Cmu$---the ``size'' of an
\eM\ measured as the Shannon entropy of the set of causal states---measures
the amount of historical information that the source stores. That
is, \eM\ minimality is not only helpful in terms of compact representations,
but it is essential, since $\Cmu$ gives one a quantitative way to say
how structured a hidden information source is.

Although the emphasis here is on the mathematical foundations of computational
mechanics, we should note that it has been used to analyze structural
complexity in a wide range of information sources. These include cellular
automata,~\cite{h93} one-dimensional maps,~\cite{Crut88a,cy90} and the
one-dimensional Ising model,~\cite{f98a,cf97} as well as several experimental
systems, such as the dripping faucet,~\cite{gpso98} atmospheric
turbulence,~\cite{pfb00} geomagnetic data,~\cite{cfw01} complex
materials,~\cite{v01,vcc02a} and molecular dynamics.~\cite{nerukh03a,nerukh03b}

In the present work we begin to address the problems posed in Appendix H.3 of
Ref. \onlinecite{Shal98a} on founding computational mechanics more fully on
stochastic process and measure theories by considering one part of the Forward
Modeling Problem noted above. The results here differ from previous work on
computational mechanics in two ways. First, the development is mathematically
rigorous, in the sense that we use measure theory to explore the notion of
minimal representations, which underlies \eMs. What is novel compared to
stochastic process theory is that we ask for minimal representations of a
stochastic process and express them in terms of the minimal \sAlg. We also
introduce two new components of the minimization procedure---internal-event
and internal-state reduction---which complement the existing concept of
causal-state reduction for \eMs.
Analyzing the Forward Modeling Problem in this way allows us to draw parallels
with the computational mechanics development of \eMs, comparing and contrasting
the various kinds of reduction method. We show that in a number of cases
these reductions are equivalent and so provide an extension of the original
concept of an \eM\ to a broader class of processes than previously possible.

\subsection{Related Work}

The modeling questions that we address here, and that are also addressed by
computational mechanics, do not arise in a vacuum. Here we briefly mention
related work that is motivated by similar concerns of the equivalence of
observed processes and of structural inference, but that adopts different
approaches. In a later section, when we turn to discuss our results, we
broaden the discussion of related work to mention additional areas in which
one might find useful applications.

One of the first attempts to address the difficulties of analyzing (known)
hidden information sources is that of Ref. \onlinecite{Blac57a}. The problem,
which comes under the heading of the \emph{identifiability of functions of
Markov chains}, was to calculate the source entropy rate, given an internal
finite-state Markov chain, the states of which are observed with a
probabilistic measurement function. (Note that today one refers to this class
of information sources as \emph{hidden Markov models}.\cite{Rabi89a,Elli95a})
There it was shown that in the majority of cases there are no closed-form
expressions for the entropy rate. A corollary of this result is that one needs
to determine the effective states (and these might be infinite in number) in
order to calculate a property as basic as the entropy rate---that is, simply
attempting to determine how random a finite-state information source is.
This contrasts, of course, with Shannon's closed-form expression for finite
Markovian sources.\cite{Shan62} We take Ref. \onlinecite{Blac57a}'s result
as one of the first indications of the nontrivial nature of inferring the
structure of hidden information sources. Another testimony to this difficulty
is that the problem of identifiability itself, though posed by Blackwell and
Koopmans in the late 1950s, was not solved for almost 40 years.\cite{Ito92a}
Moreover, the existence of minimal representations of these same hidden
sources was not established until a few years later still.\cite{Crut92c,Uppe97a}

Similar concerns about inference, representation, and causality are found in
the fields of causal inference,\cite{Glym99a} graphical models,\cite{Jord99a}
and nonlinear time series analysis and state-space reconstruction.\cite{Casd91a}
Most of the work in these areas proceeds by \emph{assuming} a given set of
observed and hidden variables (and their connectivity) and then asks for
efficient algorithms to estimate various kinds of marginal, conditional, and
joint distributions. The goals are to infer from the latter the relationship
between these variables and so, on that basis, to draw structural conclusions.
That is, in these cases one begins with strong structural priors about the
internal architecture of a hidden information source in order to initiate
analysis. Notably, only the last of these fields concentrates on temporal
dynamics and sources with memory. Here we are interested in both architectural
and temporal properties of memoryful hidden sources and wish to understand
these employing a minimum of structural priors.

\subsection{Outline}

The principle focus of the following is to develop the notion of a minimal
reduction of a given (hidden) Markov process. To do this, the development
is organized as follows. In the next section we characterize the
(rather general) class of stochastic processes---hidden information
sources---in a way that respects the distinction between a process's
internal structure and the measurements, which indirectly reflect the
internal state, available to an observer. This then allows us to define
generators of stochastic processes as \emph{Markov transition kernels},
and so state the problem of observationally equivalent generators.
The succeeding section establishes how different generators
can be mapped onto each other while maintaining observational
equivalence. Then, in the next section, we address the central problem
and show that one can maximally reduce the representation of a
process's internal structure---it's generator---while still producing
the same observed stochastic process. The reduction is achieved in two
steps---the first, called \emph{internal-event reduction}, produces
the smallest \sAlg\ and the second, \emph{internal-state reduction},
reduces the internal structure further, removing components that
are not necessary for optimal prediction. During the development we
illustrate the ideas with several examples that show how the
new formulation extends the range of applicability of computational
mechanics.

\section{Generators of Stochastic Processes}
\label{Sec:Generators}

An information source is a process that at each time step emits an output or
measurement symbol. Only the probabilistic nature of the output process is
specified in order to describe the observed information processing. Indeed,
often in information theory a source is mathematically described as a
stochastic process without concrete specification of internal mechanisms.
In many theories of complexity, however, one often uses explicitly structural
notions (e.g., \emph{automata}) from the theory of discrete
computation\cite{Hopp&Ull} to describe the resources required to reproduce or
model an observed process. So that we will have a mathematical model that
both captures the observed stochastic process and allows for a range of internal
structures, we adapt the concept of finite-state automata to the setting of
stochastic processes as follows; cf. Refs. \onlinecite{Paz71a} and
\onlinecite{Elli95a}.

\begin{figure}
\begin{center}
\setlength{\unitlength}{1cm}
\begin{picture}(6,6)
\put(0,0){\epsfxsize=6cm\epsfbox{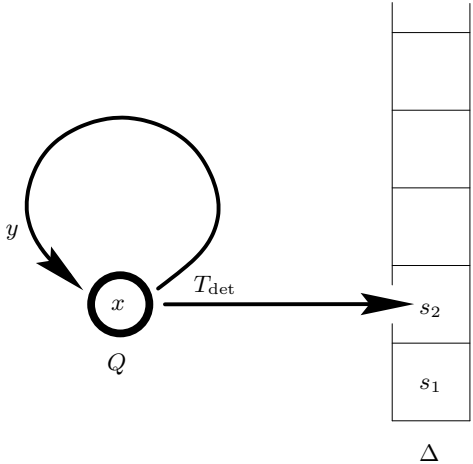}}
\put(1.15,0.7){$Q$}
\put(2.3,1.75){$T_{\rm det}$}
\put(1.2,1.5){$x$}
\put(-0.2,2.5){$y$}
\put(5.3,-0.5){$\OSym$}
\put(5.3,0.45){$s_1$}
\put(5.3,1.45){$s_2$}
\end{picture}
\end{center}   
\caption{The transition structure of a deterministic machine that
  generates a stochastic output process.
  }
\label{fig:TransitionStructure}
\end{figure}

We consider a finite set $Q$ of \emph{internal states} of the system and also
a finite set $\OSym$ of \emph{output states}, which are the observed symbols. 
The internal structure is modeled in various ways. First, it can be specified
by a deterministic (\emph{det}) transition map:
\begin{equation}
\label{det}
    T_{\rm det} : \;\; Q \; \to \; Q \times \OSym,
        \qquad x \; \mapsto \; T_{\rm d}(x) = (y,s) ~.
\end{equation}     
This map assigns to each internal state $x \in Q$ the next internal state $y$
and, at the same time, also the next output symbol $s \in \OSym$. Figure
\ref{fig:TransitionStructure} illustrates the transition structure.

A nondeterministic (\emph{non}) version of such a
\emph{machine (without input)} can be introduced as a map
\begin{equation} \label{nondet}
  T_{\rm non}: \;\; Q \; \to \; 2^{Q \times \OSym} ~,
  \qquad x \mapsto C ~.
\end{equation}
This machine assigns to each internal state $x$ a set $C$ of possible
next state pairs $(y,s)$. This extends the deterministic machine
$T_{det}$ of Eq. (\ref{det}), which can be interpreted within the
nondeterministic framework as follows:
\begin{equation*}
  T_{\rm non}(x) \; := \;
  \{T_{\rm det}(x)\} \; \in \; 2^{Q \times \OSym} ~. 
\end{equation*}

Finally, a further extension is provided by the following probabilistic
(\emph{pr}) interpretation of a nondeterministic machine $T_{\rm non}$: 
\begin{equation*}
  T_{\rm pr}: \;\; Q \times 2^{Q \times \OSym} \; \to \; [0,1] ~,
\end{equation*}
where
\begin{equation*}
  (x,C) \; \mapsto \; T_{\rm pr}(x,C) :=
  \frac{|C \cap T_{\rm non}(x)|}{|T_{\rm non}(x)|} ~.
\end{equation*}
The function $T_{\rm pr}$ satisfies  
\begin{equation*}
  T_{\rm pr}(x,C_1 \uplus C_2) \; = \;
  T_{\rm pr}(x,C_1) + T_{\rm pr}(x, C_2)
\end{equation*}
and
\begin{equation*}
  T_{\rm pr}(x,Q \times \OSym) \; = \; 1 ~.
\end{equation*}
Therefore $T_{\rm pr}$ is a \emph{Markov transition kernel} on finite symbols.

This interpretation allows for an extension of finite-state machines to
machines given by general Markov transition kernels that are not restricted
to finite symbols, for example. Here, we
allow the internal states to be described by an arbitrary measurable
space $(Q,{\cal Q})$. Again, $Q$ is the set of internal states or, in
terms of probability theory, the set of (internal) elementary events.
The $\sigma$-algebra ${\cal Q}$ represents all internal events of
interest.  The output is modeled by a measurable space
$(\OSym,{\cal D})$, too. A machine is now considered to be a Markov
transition kernel:
\begin{equation*}
   T: \;\; Q \times \big({\cal Q} \otimes {\cal D}\big) ~,
   \qquad (x,C) \; \mapsto \; T(x,C) ~.
\end{equation*}
More precisely, $T$ is assumed to satisfy the following conditions:
\begin{enumerate}
\item For all $x \in Q$, the function $T(x,\cdot)$ is a probability
        distribution on ${\cal Q} \otimes {\cal D}$.
\item For all $C \in {\cal Q} \otimes {\cal D}$, the function
        $T(\cdot,C)$ is ${\cal Q}$-measurable.
\end{enumerate}

We should point out that the well established notion of
machines that manipulate finitely many (or a countable number of)
symbols may seem more appropriate for implementations in physical
systems than our broad approach to computation using general Markov
transition kernels. Putting the natural ideas of computation into the 
probabilistic setting, however, allows us to employ measure-theoretic
concepts and techniques. This approach turns out to be very useful
in understanding the relations between the probabilistic nature of
the observed processes and the underlying internal computational
structures processes. In particular, problems on minimality properties
of machines can be handled in an efficient way and for a broader
class of processes than those over discrete symbols.

\begin{figure}
\begin{center}
\setlength{\unitlength}{1cm}
\begin{picture}(5,5)
\put(0,0){\epsfxsize=5cm\epsfbox{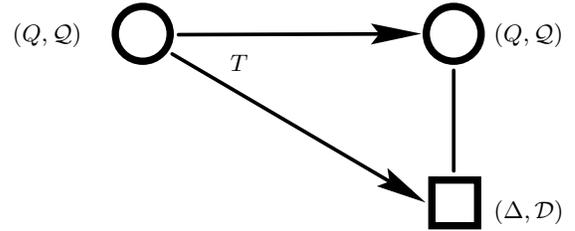}}
\put(-1.3,2.55){$(Q,{\cal Q})$}
\put(5.1,2.55){$(Q,{\cal Q})$}
\put(5.1,0.2){$(\OSym,{\cal D})$}
\put(1.6,2.15){$T$}
\end{picture}
\end{center}   
\caption{The Markov transition kernel: The internal structure is
  specified by $(Q,{\cal Q})$ and the observed process by
  $(\OSym,{\cal D})$.
  }
\label{fig:MarkovTransitionKernel}
\end{figure}

\begin{figure}
\begin{center}
\setlength{\unitlength}{1cm}
\begin{picture}(8,4)
\put(0,1){\epsfxsize=8cm\epsfbox{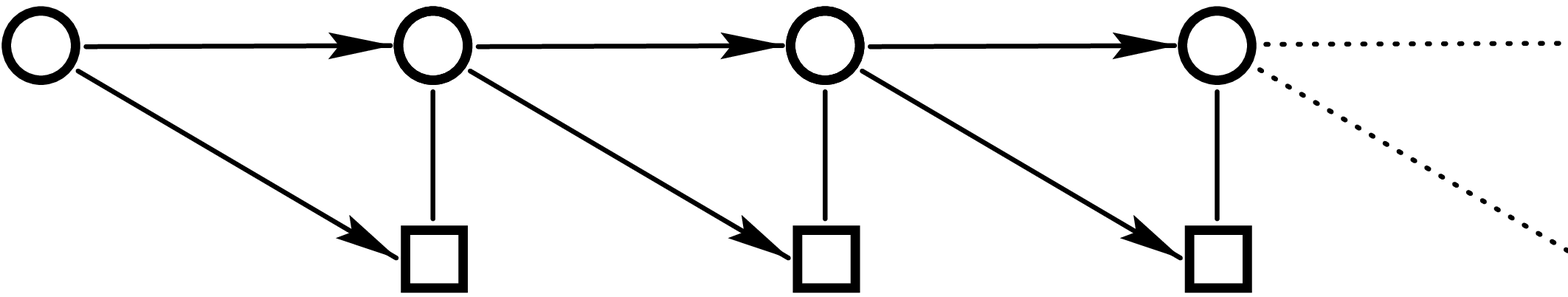}}
\put(0.125,2.6){$t = 0$}
\put(0.1,1.8){$\mu$}
\put(2.1,2.6){$t = 1$}
\put(2,0.65){$B_1$}
\put(4.1,2.6){$t = 2$}
\put(4,0.65){$B_2$}
\put(6.05,2.6){$t = 3$}
\put(5.97,0.65){$B_3$}
\end{picture}
\end{center}   
\caption{The finite-dimensional marginals ${\rm P}_n^{\mu,T}$ on
  $(\OSym^n, {\cal D}^{n})$.
  }
\label{fig:FiniteDimensionalMarginals}
\end{figure}

Given a Markov transition kernel $T$ from $(Q,{\cal Q})$ to
$(Q\times \OSym, {\cal Q}\otimes {\cal D})$, we consider it as a
temporal ``map'', as illustrated in Fig.
\ref{fig:MarkovTransitionKernel}.
In order to specify observable stochastic processes in
$(\OSym,{\cal D})$, we consider an initial distribution $\mu$ on
$(Q,{\cal Q})$ and measurable sets $B_1,\dots,B_n \in {\cal D}$. 
The finite-dimensional marginals
${\rm P}_n^{\mu,T}$ on $(\OSym^n, {\cal D}^{n})$ are obtained by
iteration of $T$, as shown in Fig. \ref{fig:FiniteDimensionalMarginals}.

This suggests the following expression for the finite-dimensional
marginals of the observed stochastic process:
\begin{widetext}
\begin{equation} \label{visible}
{\rm P}_n^{\mu, T} (B_1 \times \cdots \times B_n) \\ \nonumber
   := \int_{Q}\int_{Q \times B_1}\cdots\int_{Q\times B_n} 
   T(x_{n-1}, d(x_n,y_n)) \cdots T(x_0 , d(x_1,y_1)) \, \mu(d x_0) ~. 
\end{equation}
\end{widetext}
(Throughout the following $d(x,y)$ denotes the differential of two
variables. This notation should not be confused with a distance
measure between $x$ and $y$.)\\

\noindent
{\bf Proposition 2.1.} {\em Up to equivalence, there is exactly one
stochastic process $Y_n$, $n=1,2,\dots$, in $(\OSym,{\cal D})$, such
that for all
$n \in {\Bbb N}$ and all $B_i \in {\cal D}$, $i = 1, \dots, n$, 
\begin{equation} \label{finite}
  {\rm Pr}\{Y_1 \in B_1, \dots, Y_n \in B_n\} \;
        = \; {\rm P}_n^{\mu,T}(B_1 \times \cdots \times B_n) ~.
\end{equation}
We can identify this process, or more precisely, the class of
corresponding equivalent processes, with a probability distribution
${\rm P}^{\mu, T}$ on $(\OSym^{\Bbb N},{\cal D}^{\Bbb N})$. 
}\\
\noindent
{\em Proof.} This follows from Kolmogorov's extension
theorem.\cite{Baue72a} \hfill $\Box$\\

\noindent
{\bf Definition 2.2.} We call a Markov transition kernel $T$ from 
$(Q,{\cal Q})$ to $(Q \times \OSym, {\cal Q} \otimes {\cal D})$ a
\emph{generator} and denote it by
$[(Q,{\cal Q}),T,(\OSym,{\cal D})]$ or simply by $T$. 
We say that a stochastic process $(Y_n)_{n \in {\Bbb N}}$ in
$(\OSym,{\cal D})$ is \emph{generated} by $T$ if there exists a
probability distribution $\mu$ on $(Q,{\cal Q})$ such that  
Eq. (\ref{finite}) is satisfied for all 
$n \in {\Bbb N}$ and all $B_i \in {\cal D}$, $i=1,\dots,n$. 

Given a stochastic process $Y = (Y_n)_{n \in {\Bbb N}}$, a natural
question is whether there always exists a generator that generates $Y$.
The following trivial shift ansatz shows that this is indeed the
case.\\

\noindent
{\bf Example 2.3 (Shift Generator):}
We set $Q := \OSym^{\Bbb N}$ and ${\cal Q} := {\cal D}^{\Bbb N}$.
Consider the shift map $s: Q \to Q$, where
\begin{equation*}
  x = (y_n)_{n \in {\Bbb N}} \;
  \mapsto  \; s(x) = (y_{n+1})_{n \in {\Bbb N}} ~,
\end{equation*}
and the projection onto the first coordinate $\pi: Q \to \OSym$,
where
\begin{equation*}
  x = (y_n)_{n \in {\Bbb N}} \; \mapsto  \; \pi(x) = y_1 ~.
\end{equation*}
Furthermore, we define
\begin{eqnarray*}
  T(x, A \times B) 
  & := & \left\{
          \begin{array}{c@{,\quad}l}
            1 & \mbox{if $s(x) \in A$ and $\pi\big(s(x)\big)\in B$} \\
            0 & \mbox{otherwise} ~.
          \end{array}
        \right. 
\end{eqnarray*}
Now, a stochastic process $Y = (Y_n)_{n\in {\Bbb N}}$ in 
$(\OSym,{\cal D})$ can be identified with a probability distribution
$\mu$ on $(Q,{\cal Q}) = (\OSym^{\Bbb N}, {\cal D}^{\Bbb N})$. It is
easy to prove that $T$ generates $Y$ by verifying Eq. (\ref{finite})
with initial distribution $\mu$.\hfill $\Box$

The shift generator (Example 2.3) is maximal in the
sense that it generates
all processes in $(\OSym,{\cal D})$. For an arbitrary generator $T$,
we consider the map
\begin{equation*}
  G_T: \;\; \boldsymbol{P}(Q,{\cal Q}) \;
  \to \; \boldsymbol{P}(\OSym^{\Bbb N},{\cal D}^{\Bbb N}), \qquad 
   \mu \; \mapsto \; G_T(\mu) ~.
\end{equation*}
(Throughout, for a general measurable space $(X,{\cal X})$,
$\boldsymbol{P}(X,{\cal X})$ denotes the set of probability measures
on $(X,{\cal X})$.) The image ${\rm im}(G_T)$ of $G_T$ is the set of
processes that are generated by $T$. Here, we mainly focus on the
following problem.\\

\noindent
{\bf Problem Statement 2.4:}
{\em Given a generator $T$, can we find
a substitute $T'$ for $T$, which, on the one hand, generates the same 
set of processes, that is ${\rm im}(G_T) = {\rm im}(G_{T'})$, and, on
the other, is minimal in some sense?}\\

From Eq. (\ref{visible}) it follows directly that $G_T$ is affine in
the sense that for all $\mu_1,\mu_2 \in \boldsymbol{P}(Q,{\cal Q})$
and all $0\leq t \leq 1$, 
\begin{equation} \label{generierung}
  G_T\big((1-t) \, \mu_1 + t \, \mu_2\big) \;
  = \; (1-t) \, G_T (\mu_1) + t\, G_T(\mu_2) ~. 
\end{equation}
This implies that ${\rm im}(G_T)$ is a convex set, and we have the
following constraint on the solution of Problem 2.4:
The set
${\rm ext}\big({\rm im}(G_T)\big)$ of the extreme points of
${\rm im}(G_T)$ represents a ``lower bound'' for the set $Q$ of
internal states. 
More precisely, we have the following onto mapping $Q \to {\rm ext}\big({\rm im}(G_T)\big)$:
\begin{equation*}
  x \;\; \mapsto \;\; \delta_x \;\; \mapsto \;\; G_T(\delta_x) ~.
\end{equation*}
Thus, we cannot expect to have a notion of minimality that reduces the
internal states more than given by the extreme points of
${\rm im}(G_T)$.

However, identifying internal states $x_1$ and $x_2$
if $G_T(\delta_{x_1}) = G_T(\delta_{x_2})$ leads to a partition of $Q$
into equivalence classes---classes that are the analogs of the
\emph{causal states} in computational mechanics. The corresponding
canonical projection of internal states to their equivalence classes
is called \emph{causal-state reduction}, which is intended to reduce
the internal structure in such a way that a given observed stochastic
process is still generated by the reduced generator.

This is different from the intention stated in Problem 2.4, which is
to reduce a given generator without affecting the \emph{whole} set of
observable stochastic processes. We solve
this problem by applying reductions within a natural category of
generators. The morphisms of this category will be introduced in
Section \ref{Transformation}. Based on the results there, we present
our reduction procedures in Section \ref{Reductions}. We leave to
the future discussing causal-state reduction in terms of
morphisms in a larger category than the one studied here.
 
\section{Transformation Rules for Generators}
\label{Transformation}

\begin{figure}
\setlength{\unitlength}{1cm}
\begin{center}
\begin{picture}(6,10)
\put(0,1){\epsfxsize=6cm\epsfbox{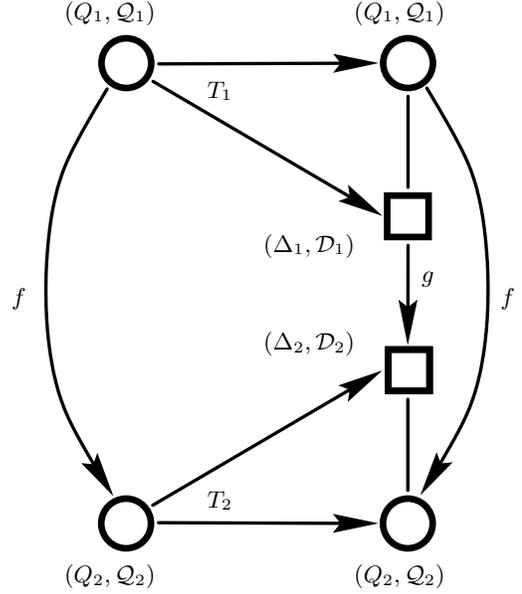}}
\put(0.3,8.1){$(Q_1,{\cal Q}_1)$}
\put(4.15,8.1){$(Q_1,{\cal Q}_1)$}
\put(2.2,7.05){$T_1$}
\put(0.3,0.6){$(Q_2,{\cal Q}_2)$}
\put(4.15,0.6){$(Q_2,{\cal Q}_2)$}
\put(2.2,1.6){$T_2$}
\put(2.95,5){$(\OSym_1,{\cal D}_1)$}
\put(2.95,3.7){$(\OSym_2,{\cal D}_2)$}
\put(-0.4,4.3){$f$}
\put(6.1,4.3){$f$}
\put(5.05,4.6){$g$}
\end{picture}
\end{center}
\caption{Commutativity for generators of equivalent observed processes.
  }
\label{fig:Commutativity}
\end{figure}

We interpret generators as objects of a category and define the
morphisms between these objects in the following way: Let
$[(Q_i,{\cal Q}_i),T_i,(\OSym_i,{\cal D}_i)]$, $i=1,2$, be two
generators. A morphism $T_1 \to T_2$ consists of  a pair $(f,g)$ of
measurable maps $f: Q_1 \to Q_2$ and $g: \OSym_1 \to \OSym_2$ such 
that for all $x \in Q_1$, $A \in {\cal Q}_2$, and $B \in {\cal D}_2$
the following commutativity rule holds:
\begin{equation} \label{consistenz}
  T_2(f(x), A \times B) \;
  = \; T_1\left(x, f^{-1}(A) \times g^{-1}(B)\right) ~.
\end{equation}
The diagram in Fig. \ref{fig:Commutativity} illustrates this
commutativity.

With the product map
\begin{equation*}
  (f \times g):\;\; Q_1 \times \OSym_1 \;
  \rightarrow \;  Q_2 \times \OSym_2, \qquad (x,y) \; 
    \mapsto \; (f(x),g(y)) ~, 
\end{equation*}
we can rewrite Eq. (\ref{consistenz}) as
\begin{equation*}
  T_2(f(x), A \times B) \;
  = \; T_1\big(x, (f \times g)^{-1}(A \times B)\big) ~.
\end{equation*}
Thus, the property of Eq. (\ref{consistenz}) is equivalent to
\begin{equation} \label{andere1}
   T_2(f(x), C) \; = \; T_1\big(x,(f \times g)^{-1}(C)\big) ~,
\end{equation}
with $C \in {\cal Q}_2 \otimes {\cal D}_2$.
Here, one has to use the fact that two probability measures are equal
if they coincide on an intersection closed system of measurable sets
that generates the underlying $\sigma$-algebra.\cite{Baue72a}
Rewriting (\ref{andere1}) gives us 
\begin{equation} \label{andere}
   T_2(f(x), \cdot) \; = \; (f\times g)_\ast \big(T_1(x,\cdot)\big) ~,
\end{equation}
where $(f\times g)_\ast(\mu)$ denotes the $(f\times g)$-image of a
probability distribution $\mu$.

In the following, a morphism $(f,g)$ is called a
\emph{transition-preserving map}. In order to define the composition
of transition-preserving maps, we consider three generators
$[(Q_i,{\cal Q}_i),T_i,(\OSym_i,{\cal D}_i)], ~ i = 1,2,3$,
and transition-preserving maps $(f_i,g_i): T_i \to T_{i+1}$, $i=1,2$. 
Now define the composition as
\[
(f_2,g_2) \circ (f_1,g_1) \; := \; (f_2\circ f_1, g_2\circ g_1) ~. 
\]
We prove that this composition is a transition-preserving map
$T_1 \to T_3$ by verifying Eq. (\ref{consistenz}):
\begin{eqnarray*}
T_3\big((f_2 & \circ & f_1)(x), A \times B \big)
= T_3\Big( f_2\big(f_1(x)\big),A \times B\Big) \\
& = & T_2\big( f_1(x), f_2^{-1}(A) \times g_2^{-1}(B)\big) \\
& = & T_1\Big( x, f_1^{-1}\big(f_2^{-1}(A)\big) \times g_1^{-1}\big(g_2^{-1}(B)\big) \Big) \\
& = & T_1\big( x, (f_2 \circ f_1)^{-1}(A) \times (g_2 \circ g_1)^{-1}(B)\big) ~.
\end{eqnarray*}

\noindent
{\bf Proposition 3.1.} {\em Let $[(Q_i,{\cal Q}_i),T_i,(\OSym_i,{\cal D}_i)]$, $i=1,2$, 
be two generators, and let $(f,g)$ be a transition-preserving map from $T_1$ to $T_2$, 
and let $\mu$ be a probability distribution on $(Q_1,{\cal Q}_1)$. Then, denoting the $f$-image of $\mu$ by 
$f_\ast(\mu)$, for all $B_1,\dots,B_n \in {\cal D}_2$,
\[
{\rm P}_n^{f_\ast(\mu),T_2}(B_1 \times \cdots \times B_n)  
\; = \; {\rm P}_n^{\mu,T_1}(g^{-1}(B_1) \times \cdots \times g^{-1}(B_n)) ~.
\]
}  
\noindent
{\em Proof.} With the general transformation rule for integrals we have
\begin{widetext}
\begin{eqnarray*}
{\rm P}_n^{f_\ast(\mu),T_2}(B_1 \times \cdots \times B_n)
& = & \int_{Q_2}\int_{Q_2 \times B_1}\cdots\int_{Q_2\times B_n} 
      T_2(x'_{n-1},d(x'_n,y'_n)) \, \cdots \, T_2(x'_0, d(x'_1,y'_1)) \, f_\ast(\mu)(d x'_0) \\
& = & \int_{Q_1}\int_{Q_2 \times B_1}\cdots \int_{Q_2\times B_n} 
      T_2(x'_{n-1},d(x'_n,y'_n)) \, \cdots \, T_2(f(x_0), d(x'_1,y'_1)) \, \mu(d x_0) \\
&   & \mbox{(transformation rule)} \\
& \stackrel{(\ref{andere})}{=} & \int_{Q_1}\int_{Q_2 \times B_1}\cdots \int_{Q_2\times B_n} 
      T_2(x'_{n-1},d(x'_n,y'_n)) \,  
      \cdots \, (f \times g)_\ast\big(T_1(x_0, \cdot)\big) (d(x'_1,y'_1)) \, \mu(d x_0) \\
& = & \int_{Q_1}\int_{Q_1 \times g^{-1}(B_1)}\cdots \int_{Q_2\times B_n} 
      T_2(x'_{n-1},d(x'_n,y'_n)) \, \cdots \, T_1(x_0, d(x_1,y_1)) \, \mu(d x_0) \\
&   & \mbox{(transformation rule)} \\
& \vdots & \qquad \vdots \qquad \qquad \vdots \\
& = & \int_{Q_1}\int_{Q_1 \times g^{-1}(B_1)}\cdots \int_{Q_2\times g^{-1}(B_n)} 
      T_1(x_{n-1},d(x_n,y_n)) \, \cdots \, T_1(x_0, d(x_1,y_1)) \, \mu(d x_0) \\
& = & {\rm P}_n^{\mu,T_1}(g^{-1}(B_1) \times \cdots
   \times g^{-1}(B_n)) ~.
\end{eqnarray*}
$ $ \hfill $\Box$
\end{widetext}

\noindent
{\bf Theorem 3.2.} {\em If $T_1$ generates $(Y_n)_{n \in {\Bbb N}}$
and $(f,g)$ is a transition-preserving map from $T_1$ to $T_2$, then 
$T_2$ generates $(g \circ Y_n)_{n \in {\Bbb N}}$.
}\\

\noindent
{\em Proof.} This statement follows directly from Proposition 3.1.
\hfill $\Box$

Theorem 3.2 has important and direct implications for two special cases.
In the first case, we fix $g$ as the identity map and, in the second
case, we fix $f$ as the identity map. In these cases, without reference
to the identity maps, $f$ and $g$ are called transition-preserving.
The implications are stated in the following two corollaries.\\

\noindent
{\bf Corollary 3.3.} {\em Let $[(Q_i,{\cal Q}_i),T_i,(\OSym,{\cal D})]$,
$i=1,2$, be two generators, and let $f$ be a transition-preserving map.
Then
\[
   G_{T_1} \; = \; G_{T_2} \circ f_\ast ~.
\]
In particular, this implies 
\[
  {\rm im}(G_{T_1}) \; \subseteq \; {\rm im}(G_{T_2}) ~,  
\]
where the equality holds if $f_\ast$ is onto.}\\

\noindent
{\bf Corollary 3.4.} {\em 
Let $[(Q,{\cal Q}),T_1,(\OSym_1,{\cal D}_1)]$ be a generator of a
stochastic process $(Y_n)_{n\in {\Bbb N}}$ in $(\OSym_1,{\cal D}_1)$, 
and let $g: (\OSym_1,{\cal D}_1) \to (\OSym_2,{\cal D}_2)$ be a
measurable map. Then
$T_2: Q \times ({\cal Q} \otimes {\cal D}_2) \to [0,1]$ with 
\[
   T_2(x,A \times B) \; := \; T_1\big(x, A \times g^{-1}(B)\big)
\]
is a generator of the stochastic process
$(g \circ Y_n)_{n \in {\Bbb N}}$ in $(\OSym_2,{\cal D}_2)$.
}

\section{Reductions of Generators}
\label{Reductions}

After having derived some basic transformation rules for generators in
Section \ref{Transformation}, we are now ready to concentrate on the
main problem, namely to maximally reduce a given generator $T$ while
keeping the set ${\rm im}(G_T)$ of generated processes unchanged. The
solution of this problem is given by Theorem 4.5 below and is based on a
combination of reduction methods, which we present in this section.
First, we attempt to reduce the $\sigma$-algebra ${\cal Q}$ of internal
events as much as possible, by considering only those events in
${\cal Q}$ that are necessary for maintaining the output process
unchanged. The following theorem formalizes this idea.\\

\noindent
{\bf Theorem (Internal-Event Reduction) 4.1.} {\em Let 
$[(Q,{\cal Q}),T, (\OSym,{\cal D})]$ be a generator. Then there 
exists a smallest $\sigma$-subalgebra $\sigma_Q(T)$ of ${\cal Q}$ with the property 
that for all $C \in \sigma_Q(T) \otimes {\cal D}$, $T(\cdot, C)$ is 
$\sigma_Q(T)$-measurable. 
The generator $[(Q,\sigma_Q(T)),\bar{T}, (\OSym,{\cal D})]$
with the restriction 
$\bar{T} := T|_{Q \times (\sigma_Q(T) \otimes {\cal D})}$ then satisfies 
\[
     {\rm im}(G_{\bar{T}}) \; = \; {\rm im}(G_T) ~.
\]
}
\noindent
{\em Proof.} Let ${\cal A}_i$, $i \in I$, be the family of all $\sigma$-subalgebras of ${\cal Q}$ that satisfy
the following condition: 
for all $C \in {\cal A}_i \otimes {\cal D}$, $T(\cdot, C)$ is 
${\cal A}_i$-measurable. Now define
\[
   \sigma_Q(T) \; := \; \bigcap_{i \in I} {\cal A}_i ~.  
\]  
Then for $C \in \sigma_Q(T) \otimes {\cal D}$, $T(\cdot,C)$ is ${\cal A}_i$ measurable for 
all $i \in I$, and therefore also $\sigma_Q(T)$-measurable.

For the reason that trivially
\[
   \bar{T}({\rm id}_Q(x), A \times D)
   \; = \; T(x, {\rm id}_Q^{-1}(A) \times {\rm id}_\OSym^{-1}(D)) ~, 
\]
Corollary 3.3 implies that $T$ and $\bar{T}$ generate the same set of
stochastic processes. \hfill $\Box$\\

\noindent
Theorem 4.1 guarantees the existence of a minimal sufficient $\sigma$-subalgebra of ${\cal Q}$. 
Now we provide a way to calculate it explicitly in the case where we have a deterministic internal dynamics 
$f: Q \to Q$ and a visible process given by a measurement $g: Q \to \OSym$. This case generalizes the shift generator 
of Example 2.3.    
\vspace{1cm}

\noindent
{\bf Theorem 4.2.} {\em Let $(Q,{\cal Q})$ and $(\OSym,{\cal D})$ be two measurable spaces, and let 
$f:Q \to Q$ and $g:Q \to \OSym$ be two measurable maps. Consider the generator $[(Q,{\cal Q}),T,(\OSym,{\cal D})]$ 
defined by
\begin{eqnarray*}
  T(x,A \times B) \; & := & \; 1_{f \in A,\, g \circ f \in B}(x) \\
  \; & = & \;
    \left\{
      \begin{array}{c@{,\quad}l}
        1 & \mbox{if $f(x) \in A$ and $g\big(f(x)\big) \in B$} \\
        0 & \mbox{otherwise}
      \end{array}
    \right. ~.
\end{eqnarray*}
Then
\begin{equation} \label{explizitsp}
       \sigma_Q(T) \; = \; \sigma(g \circ f, g\circ f^2,\dots) ~.
\end{equation}
}
\noindent
{\em Proof.} We prove inclusion in each direction separately.
\begin{enumerate}
\item We establish that
$\sigma_Q(T) \subseteq  \sigma(g \circ f, g\circ f^2,\dots)$ by showing that
for all 
$C \in \sigma(g \circ f, g\circ f^2,\dots) \otimes {\cal D}$, $T(\cdot,C)$
is measurable with respect to $\sigma(g \circ f, g\circ f^2,\dots)$. From 
\begin{eqnarray*}
  \sigma(g & \circ & f, g \circ f^2, \dots ) \otimes {\cal D} \\
  \; & = & \; \sigma\big((g\circ f, g\circ f^2, \dots) \times {\rm id}_\OSym\big)
\end{eqnarray*}
we know that there exists a measurable set
$C' \in {\cal D}^{\Bbb N} \otimes {\cal D}$ with
\[
    C \; = \; \big((g\circ f, g\circ f^2, \dots)
	\times {\rm id}_\OSym\big)^{-1}(C') ~.
\]
This implies
\begin{eqnarray*}
T(\cdot,C) & = & 1_{(f,g\circ f) \; \in \; C} \\
           & = & 1_{(f,g\circ f) \; \in \; \big((g\circ f, g\circ f^2, \dots) \times {\rm id}_\OSym\big)^{-1}(C')} \\
           & = & 1_{\big((g\circ f, g\circ f^2, \dots) \times {\rm id}_\OSym\big)\circ (f,g\circ f) \; \in \; C'} \\
           & = & 1_{((g\circ f^2,g\circ f^3,\dots), g \circ f) \; \in \; C'} ~.
\end{eqnarray*}
Thus,
$T(\cdot,C)$ is measurable with respect to $(g\circ f, g\circ f^2, \dots)$.
\item Now we prove that
$\sigma_Q(T) \; \supseteq \; \sigma(g \circ f, g\circ f^2, \dots)$
by applying an induction argument to show that $\sigma(g \circ f^k) \subseteq \sigma_Q(T)$ for all $k= 1,2,\dots$:
\begin{enumerate}
\item ``$k=1$'': Let $A$ be a $(g \circ f)$-measurable set. Then there exists
        a measurable set $B \in {\cal D}$ with $A = (g \circ f)^{-1}(B)$. From
        $Q\times B \in \sigma_Q(T)$, and
\begin{eqnarray*}
   1_A \; & = & \; 1_{g \circ f \in B} \\
          & = & \; 1_{(f,g\circ f) \in Q \times B} \\
                  & = & \; T(\cdot,Q \times B) ~,
\end{eqnarray*}
        it follows that $A$ is $\sigma_Q(T)$-measurable.
\item ``$k \to k+1$'': We assume that $\sigma( g \circ f^k)$ is a
        $\sigma$-subalgebra of $\sigma_Q(T)$, and we have to show that this is
        also true for $\sigma( g \circ f^{k+1})$. To this end, we choose a
        measurable set $A \in \sigma(g \circ f^{k+1})$. There exists a measurable
        set $B \in {\cal D}$ with $A = (g \circ f^{k+1})^{-1}(B)$, and we have
\begin{eqnarray*}
   1_A & = & 1_{g\circ f^{k+1} \in B} \\
       & = & 1_{f \in (g\circ f^k)^{-1}(B)} \\
           & = & T(\cdot, (g \circ f^k)^{-1}(B) \times \OSym ) ~.  
\end{eqnarray*}
        This implies $A \in \sigma_Q(T)$, because according to the induction
        hypothesis $(g \circ f^k)^{-1}(B) \in \sigma_Q(T)$.
\end{enumerate}
\end{enumerate}
\hfill $\Box$

\noindent
{\bf Examples 4.3.}
\begin{enumerate}
\item {\bf Complete Randomness.} Consider a probability space
$(Q,{\cal Q},\mu)$. This defines the following generator
$[(Q,{\cal Q}),T,(Q,{\cal Q})]$ which is completely random in the sense 
that the next internal state, which coincides with the next output
state, is independent of the current internal state: 
\[
  T: \;\; Q \times ({\cal Q}\otimes {\cal Q}) \to [0,1] ~,
\]
and
\[
  T(x, A\times B) \; := \; \mu(A\cap B) ~.
\]
In this case
\[
    \sigma_Q(T) \; = \; \{\emptyset, Q\}.
\]
In other words, as expected, the process has no memory. Only a single internal
event is required to generate the process $\mu^{\otimes {\Bbb N}}$ and
$\mu^{\otimes {\Bbb N}}$ is the only process in ${\rm im}(G_T)$.
\item {\bf Rotation of the Unit Circle.} Consider the unit circle 
$K = \{ x\in {\Bbb C}  :  |x| = 1\}$ and its upper half 
$A_1 = \{e^{i \, \varphi}  :  \varphi \in [0,\pi) \}$ and its lower
half $A_2 = \{e^{i \, \varphi}  :  \varphi \in [\pi, 2\, \pi) \}$.
With a number $a \in K$, we construct the generator $T$ according to
Theorem 4.2 using $f(x) = a\,x$ and $g(x) = k$ for $x \in A_k$. There
are two qualitatively different cases:
\begin{enumerate}
\item Assume that $a$ is a root of unity. Then there is a natural
number $p \not= 0$ with $a^p = 1$. This implies $f^p = {\rm id}_{K}$
and, therefore,
\begin{eqnarray*}
  \sigma_Q(T) \; & = & \; \sigma(g\circ f,g\circ f^2,\dots) \\
  & = & \sigma(g \circ f, g\circ f^2,\dots,g \circ f^{p-1}) ~.
\end{eqnarray*}
Since $g$ has just two different values, $\sigma_Q(T)$ is finite in
this case and we have an effective internal-event reduction.
\item Assume that $a$ is not a root of unity. Then $\sigma_Q(T)$ is
the Borel algebra of the unit circle, and we have no internal-event
reduction.    
\end{enumerate}
\end{enumerate}
\hfill $\Box$

In addition to the reduction method given by Theorem 4.1, we now
consider another way to reduce the generator's internal structure.
Given a generator 
$[(Q,{\cal Q}),T,(\OSym,{\cal D})]$,
we identify each two elements
$x_1,x_2 \in Q$ if $T(x_1,\cdot) = T(x_2,\cdot)$. The equivalence
class of $x$ is denoted by $[x]$. Furthermore, we define
\begin{equation*}
 [Q] \; := \; \big\{[x] \; : \; x \in Q \big\}
\end{equation*}
and
\begin{equation*}
 [{\cal Q}] \; := \; \big\{ A' \subseteq [Q] \; : \;
 [\cdot]^{-1}(A') \in {\cal Q} \big\} ~.
\end{equation*}
The $\sigma$-algebra $[{\cal Q}]$ is just the terminal algebra of the
canonical projection $[\cdot]: x \mapsto [x]$. It is easy to see that
the following transition kernel
$[T]: [Q] \times \big([{\cal Q}] \otimes {\cal D}\big) \rightarrow [0,1]$ 
is well defined
\[
     [T]([x], A' \times B) \; := \; T(x, [\cdot]^{-1}(A') \times B) ~.
\]
$ $

\noindent
{\bf Theorem (Internal-State Reduction) 4.4.}
{\em Let $[(Q,{\cal Q}),T,(\OSym,{\cal D})]$ be a generator. Then 
$[([Q],[{\cal Q}]), [T],(\OSym,{\cal D})]$ is a generator, which
generates the same set of processes in $(\OSym,{\cal D})$ as $T$,
that is,
\[
     {\rm im}(G_{[T]}) \; = \; {\rm im}(G_{T}) ~. 
\]
}\\
\noindent
{\em Proof.} We show that $[T]$ is a Markov transition kernel in
two stages.
\begin{enumerate}
\item We fix $[x]$ and prove that $[T]\left([x],\cdot \right)$ is a
probability measure: 
\begin{eqnarray*}
[T]\left([x], \biguplus_{n=1}^\infty C_n \right) 
& = & T\left(x, \big([\cdot] \times {\rm id}_{\OSym}\big)^{-1}\left(\biguplus_{n=1}^\infty C_n\right)\right) \\
& = & T\left(x, \biguplus_{n=1}^\infty \big([\cdot] \times {\rm id}_{\OSym}\big)^{-1}(C_n)\right) \\ 
& = & \sum_{n = 1}^\infty
      T\left(x, \big([\cdot] \times {\rm id}_{\OSym}\big)^{-1}(C_n)\right) \\  
& = & \sum_{n = 1}^\infty
      [T]\big([x], C_n \big) ~, \\
\end{eqnarray*}
and
\begin{eqnarray*}
[T]\left([x], [Q] \times \OSym \right) 
& = & T\left(x, \big([\cdot]\times {\rm id}_\OSym \big)^{-1}\big([Q] \times \OSym\big) \right) \\
& = & T\left(x, Q \times \OSym \right) \\
& = & 1 ~. 
\end{eqnarray*}
\item Now we fix $C \in [{\cal Q}] \otimes {\cal D}$ and prove that
$[T](\cdot,C)$ is $[{\cal Q}]$-measurable. To this end, it is
sufficient to prove that for all $\varepsilon$ with 
$0 \leq \varepsilon \leq 1$, 
the set $\{[T](\cdot, C) \leq \varepsilon \}$ is an element of
$[{\cal Q}]$ or equivalently 
$[\cdot]^{-1}\big(\{[T](\cdot, C) \leq \varepsilon \}\big) \in {\cal Q}$.
This is shown as follows.
\begin{eqnarray*}
[\cdot]^{-1} \big(\{[T]( & \cdot & , C) \leq \varepsilon \}\big) \\
& = & [\cdot]^{-1}\big(\{ [x] \in [Q] \; : \; T'([x], C) \leq \varepsilon \}\big) \\
& = & \{ x \in Q \; : \; [T]([x], C) \leq \varepsilon \} \\
& = & \left\{ x \in Q \; : \; T\left(x, \big([\cdot]\times {\rm id}_\OSym\big)^{-1} (C)\right) 
      \leq \varepsilon \right\} \\
& \in & {\cal Q} ~.
\end{eqnarray*}
\end{enumerate}
\hfill  $\Box$

Combining the reduction methods provided by Theorem 4.1 and Theorem 4.4,
we can reduce every generator to a minimal generator. This statement is
specified in the following theorem.\\

\noindent
{\bf Theorem (Solution of Problem 2.4) 4.5.}
{\em Let $[(Q,{\cal Q}), T, (\OSym,{\cal D})]$ be a generator, and let 
$[ (Q',{\cal Q}'), T', (\OSym,{\cal D})]$ be the generator obtained
from $T$ by applying first the reduction method of Theorem 4.1 and
then the method of Theorem 4.4. Then $T'$ satisfies  
\[
       {\rm im}(G_{T'}) \; = \; {\rm im}(G_T)
\] 
and is minimal in the sense that given another generator
$[(Q'',{\cal Q}''),T'', (\OSym, {\cal D})]$ with 
${\rm im}(G_{T''}) = {\rm im}(G_{T'})$, every transition-preserving
map $f$ from $T'$ to $T''$ is injective.
}
\noindent
{\em Proof.} Again there are two steps.
\begin{enumerate}
\item We prove
\begin{equation} \label{gleichheit}
   \sigma(f \circ [\cdot]) \; = \; \sigma_Q(T) ~.
\end{equation}
\begin{enumerate}
\item ``$\subseteq$'': This inclusion follows directly from the
measurability of
\begin{eqnarray*}
(Q,{\cal Q}) & \stackrel{{\rm id}_Q}{\longrightarrow} & (Q,\sigma_Q(T)) \\
& \stackrel{[\cdot]}{\longrightarrow} & ([Q],[\sigma_Q(T)]) = (Q',{\cal Q}') \\
   & \stackrel{f}{\longrightarrow} & (Q'',{\cal Q}'') ~.   
\end{eqnarray*}
\item ``$\supseteq$'': Let $C \in \sigma(f \circ [\cdot])\otimes {\cal D}$. We prove that $T(\cdot,C)$ is 
$(f\circ [\cdot])$-measurable, from which $\sigma_Q(T) \subseteq \sigma(f \circ [\cdot])$ follows, because 
$\sigma_Q(T)$ is the smallest $\sigma$-algebra with that invariance property:
From
\begin{eqnarray*}
   \sigma(f \circ [\cdot])\otimes {\cal D} & = & \sigma(f \circ
   [\cdot]) \otimes \sigma({\rm id}_\OSym) \\
   & = & \sigma\big((f\circ [\cdot]) \times {\rm id}_\OSym\big) ~, 
\end{eqnarray*}
it follows that there exists $C'' \in {\cal Q}''\otimes {\cal D}$ with 
\[
  \big((f \circ [\cdot]) \times {\rm id}_\OSym\big)^{-1}(C'') \; = \; C.
\]
This implies the $(f\circ [\cdot])$-measurability of $T(\cdot,C)$:
\begin{eqnarray*}
T(x,C) 
& = & T\Big(x, \big((f \circ [\cdot]) \times {\rm id}_\OSym\big)^{-1}(C'')\Big) \\
& = & T''\big((f \circ [\cdot])(x), C''\big) \\
& = &  \big(T''(\cdot, C'') \circ  (f \circ [\cdot])\big)(x) ~.
\end{eqnarray*}
\end{enumerate}
\item Using Eq. (\ref{gleichheit}), we now prove that $f$ is injective.
Assume $f([x_1]) = f([x_2])$ where $[x_i]$ are equivalence classes in
$Q$; that is, $[x_1],[x_2] \in [Q]$. In order to prove injectivity of
$f$, we have to show $[x_1] = [x_2]$:
\begin{eqnarray*}
\bar{T}(x_1,C) 
& = & \bar{T}\Big(x_1, \big((f \circ [\cdot]) \times {\rm id}_\OSym\big)^{-1}(C'') \Big) \\
& = & T''(f([x_1]), C'') \\
& = & T''(f([x_2]), C'') \\
& = & \bar{T}\Big(x_2, \big((f \circ [\cdot]) \times {\rm id}_\OSym\big)^{-1}(C'') \Big) \\
& = & \bar{T}(x_2,C) ~. 
\end{eqnarray*}  
\end{enumerate}
$ $ \hfill $\Box$ \\

\noindent
{\bf Examples (Continuation of Examples 4.3) 4.6.}
\begin{enumerate}
\item {\bf Complete Randomness.} Applying the internal-state reduction
leads to an internal state space $Q'$ consisting of one point, namely
$Q'= \{Q\}$. The reduced generator is then given by
\[
      T'(x, \{Q\} \times B) \; = \; \mu(B) ~.  
\]
\item {\bf Rotation of the Unit Circle.}
\begin{enumerate}
\item Identifying points according to the internal-state reduction
        leads to the grouping of all elements in a given atom of the finite
        $\sigma$-algebra $\sigma_Q(T)$. Thus, in this case we have finite
        transition kernel $T'$ resulting from Theorem 4.5.
\item In this case, the internal-state reduction leads to equivalence
        classes that consist of individual points, so that effectively there
        is no reduction. 
\end{enumerate}
\end{enumerate}
$ $ \\

As pointed out at the end of Section \ref{Sec:Generators}, our goals differ
from those underlying causal-state reduction in computational mechanics.
Nonetheless, it is not hard to see the following close relationship: In the
situation of Theorem 4.2, identifying $x_1$ and $x_2$ if and only if
$G_T(\delta_{x_1}) = G_T(\delta_{x_2})$ is equivalent to the identification
of $x_1$ and $x_2$ if and only if $T(x_1,C) = T(x_2,C)$ for all
$C \in \sigma_Q(T)$. The first identification leads to the analogs of the
causal states in computational mechanics and the second identification is
the one used in Theorem 4.5. For completeness, we conclude this section with
the proof of this relationship.\\

\noindent
{\bf Corollary 4.7.} {\em Let $[(Q,{\cal Q}),T,(\Delta,{\cal D})]$ be a
generator as in Theorem 4.2, and let $x_1,x_2 \in Q$. Then 
\[
    G_T(\delta_{x_1}) \; = \; G_T(\delta_{x_2}) 
\]
is equivalent to 
\[
    T(x_1,C) \; = \; T(x_2,C) \quad \mbox{for all $C \in \sigma_Q(T)$} ~.
\]
} 
\noindent
{\em Proof.} 
\begin{eqnarray*}
  & & \quad G_T(\delta_{x_1}) = G_T(\delta_{x_2}) \\
  & \Leftrightarrow & \quad g\big(f^k(x_1)\big) = g\big(f^k(x_2)\big) \\
  &                 & \quad \mbox{for all $k=1,2,\dots$} \\
  & \Leftrightarrow & \quad 1_{(g\circ f,g\circ f^2,\dots) \in C'}(x_1) = 
                      1_{(g\circ f,g\circ f^2,\dots) \in C'}(x_2) \\
  &                 & \quad \mbox{for all $C' \in {\cal D}^{\Bbb N}$} \\
  & \Leftrightarrow & \quad T(x_1, C) = T(x_1,C) \\
  &                 & \quad \mbox{for all $C \in \sigma\{g \circ f, g \circ f^2, \dots \}$}\\
  & \Leftrightarrow & \quad T(x_1, C) = T(x_1,C) \\
  &                 & \quad \mbox{for all $C \in \sigma_Q(T)$}
  \qquad \mbox{(Theorem 4.2) ~.}
\end{eqnarray*}
\hfill $\Box$

\section{Discussion}

After this long development, it will be helpful to discuss more
informally what was achieved and how to interpret the results. We began
by characterizing the class of hidden information sources in a way that
respected the distinction between a source's internal structure and
its observed process. That allowed us to define generators of
stochastic processes as Markov transition kernels and to state the
problem of observationally equivalent generators. We then established
how different generators can be mapped onto each other while maintaining
equivalence of the observed stochastic process. We showed that one can
maximally reduce the representation of a source's generator under the same
constraint. The reduction was achieved in two steps: first by
internal-event reduction which produced the smallest \sAlg\ and the
second by internal-state reduction which collapsed \sAlg\ components
redundant for optimal prediction.

``Prediction'' here refers to the hidden internal state \emph{and} to the
observed state of the machine in the next time step. Within computational
mechanics, however, predictions are made for the whole future of the observed
process, which seems more natural than trying to make predictions of the
hidden states. For the class of generators that have the structure of Theorem
4.2 it turns out that both approaches are equivalent (see Corollary 4.7).
We expect this equivalence to be valid for a larger class of generators but
leave this to future investigations.

One interpretation of these results is that the seemingly intractable
nonuniqueness of inferring models of hidden information sources can be
directly addressed. There are more constraints on one's choice of
representation than one thinks, at first blush. The new reductions and
their sometimes-equivalence to \eM\ representations suggest that there
might be a preferred minimal representation of general stochastic
processes---the \eM\ or some generalization of it. Even if these minimal
models are unachievable when inferring from finite data, nevertheless,
they are the goal toward which modeling should strive. We hoped to show,
and partly illustrated this by the examples, that the new formulations
of reductions and their relationship to causal-state reduction greatly
extends the class of processes to which computational mechanics can be
applied.

\section{Application Areas}

The developments here properly lie in the domains of measure theory and
stochastic processes. However, we believe the results on reductions are
relevant to a number of areas outside of those fields. To emphasize
this, and also to suggest possible directions for future work, we shall
point out the similarities with some areas and possible applications
that would follow from the similarities. The areas considered are not,
by any means, exhaustive. The observations are intended only to be
suggestive.

Very generally, in statistical physics theories assume that a system
is Markovian.\cite{Penr70a} There is, for example, little
concern about minimal representations. One consequence of this is that
one sees an only indirect interest in calculating the structural and
information-processing properties of physical systems. Historically, as
reflected in the invention and use of order parameters, structural
aspects are what the theorist introduces at the beginning of analysis.
The difficulty that arises is that the systems of genuine interest
often produce ``order''---behaviors and structures---that is not
directly determined by the fundamental equations of motion, but only
arises over long times and large spatial scales. In these cases, one
must adopt something like the inferential stance to discovering the
emergent order, rather than assume it at the outset. All of which
is to say that applying the reductions discussed here to problems
in statistical mechanics should lead to novel and useful notions
of structure and to quantitative methods for measuring degrees of
structuredness.

In communication theory hidden information sources are called
\emph{channels}.\cite{Shan62} Overwhelmingly, the cases that are
considered and analyzed and that, more importantly, are the basis for
the central results of information theory assume channels with
\emph{no} memory.\cite{Cove91a} Here, though, in effect we addressed
channels with memory in the sense that the output symbols were not in
one-to-one relationship to the channel's internal states. Indeed, to
the extent the set of causal states is nontrivial, then one is
confronted with memoryful information sources. Looking forward, the
results on reductions should help in analyzing memoryful information
sources and in quantitatively addressing the size of encoders and
decoders under fixed channel fidelity.

\section{Conclusion}

The process of model building is sometimes characterized as equivalent
to data compression. While this might be true from a pragmatic
engineering perspective, from the scientific, one must disagree. Model
building is much more than data compression, especially to the extent
that one attempts to explain and understand hidden structures and
mechanisms. (See, for example, the discussion in the last section of
Ref. \onlinecite{Crut91b}.)

Building a good model certainly helps with compressing the original data,
since the predictable components of a process that the model captures can be
used in encoding and decoding to send only the ``random'' portions. However,
the goal of modeling in the sciences is understanding the (possibly hidden)
mechanisms and structures---elements that help explain observed phenomena
and lead to new insights about how nature organizes itself. In this, minimal
models---the theme of the present work---play a particularly important role.
Not only do small models make for more tractable analysis and manipulation,
they express how a process is structured and, in this, they allow for
improved scientific understanding.

Here we addressed the Forward Modeling Problem of maximally reducing a
given generator while keeping the observed process unchanged. Future work 
will focus on the Reverse Modeling Problem, the goal of which is to construct
a minimal generator based on a distribution of measurement sequences alone.
We envision a two-step approach. In the first, one constructs a possibly
large but sufficient generator that, in the second step, is reduced using
the results developed above. Unfortunately, the problem of ambiguity arises
at the end of this procedure. From previous work in computational mechanics,
however, we expect uniqueness of minimal generators up to isomorphism.

\begin{acknowledgments}

The authors thank D. Eric Smith for helpful discussions.
This work was supported at the Santa Fe Institute under the Networks Dynamics
Program funded by the Intel Corporation and under the Computation, Dynamics,
and Inference Program via SFI's core grants from the National Science and
MacArthur Foundations. Direct support was provided by DARPA Agreement
F30602-00-2-0583. NA was supported by a Santa Fe Institute Post-doctoral
Fellowship.

\end{acknowledgments}

\bibliography{chaos}

\end{document}